\documentclass{amsart}

\usepackage{amssymb}
\usepackage[all]{xy}
\usepackage{hyperref}

\newtheorem{thm}{Theorem}

\newtheorem{lem}[thm]{Lemma}
\newtheorem{cor}[thm]{Corollary}

\newtheorem{prop}[thm]{Proposition}

\newtheorem{conj}[thm]{Conjecture}
   
\theoremstyle{definition}
\newtheorem{defn}[thm]{Definition}

\newtheorem{say}[thm]{}

\newtheorem*{ack}{Acknowledgments}      

\newtheorem{defn-thm}[thm]{Definition--Theorem}  
\newtheorem{defn-lem}[thm]{Definition--Lemma}  

\theoremstyle{remark}

\let \cedilla =\c
\renewcommand{\c}[0]{{\mathbb C}}  

\renewcommand{\o}[0]{{\mathcal O}} 
\newcommand{\z}[0]{{\mathbb Z}}

\renewcommand{\r}[0]{{\mathbb R}}

\newcommand{\q}[0]{{\mathbb Q}}
\newcommand{\map}[0]{\dasharrow}
\newcommand{\qtq}[1]{\quad\mbox{#1}\quad}

\newcommand{\supp}[0]{\operatorname{Supp}}

\newcommand{\ex}[0]{\operatorname{Ex}}    
\newcommand{\diff}[0]{\operatorname{Diff}}

\newcommand{\chr}[0]{\operatorname{char}}

\newcommand{\cl}[0]{\operatorname{Cl}}

\newcommand{\rdown}[1]{\lfloor{#1}\rfloor}

\newcommand{\simr}[0]{\sim_{\r}}

\newcommand{\tsum}[0]{\textstyle{\sum}}

\def\loccoh#1.#2.#3.#4.{H^{#1}_{#2}(#3,#4)}

\DeclareMathAlphabet{\mathchanc}{OT1}{pzc}%
                                {m}{it}

\usepackage[all]{xy}\xyoption{dvips}

\begin{document}
\bibliographystyle{amsalpha}


 \title{Relative MMP without $\q$-factoriality}
\author{J\'anos Koll\'ar}

\begin{abstract} We consider the minimal model program for varieties that are not $\q$-factorial. We show that, in many cases, its steps are simpler than expected.  In particular, all flips are 1-complemented.  The main applications are to  log terminal singularities, removing the earlier $\q$-factoriality assumption from several theorems of Hacon--Witaszek   and de~Fernex--Koll\'ar--Xu. 
\end{abstract} 

 \maketitle

Let $(X, \Theta)$ be a dlt pair, projective over a base scheme $S$, and
$H$ an $\r$-divisor that is   ample over $S$. As we run the $(X, \Theta)$-MMP over $S$ 
 with scaling of $H$ as in Definition~\ref{mmp.steps.defn}, at the $i$th step 
there are 3 possibilities.
\begin{itemize}
\item (Divisorial) $X^i\stackrel{\phi_i}{\longrightarrow} Z^i= X^{i+1}$,
\item (Flipping)   $X^i\stackrel{\phi_i}{\longrightarrow}Z^i  \stackrel{\phi_i^+}{\longleftarrow}(X^i)^+=X^{i+1}$.
\item  (Mixed) $X^i\stackrel{\phi_i}{\longrightarrow} Z^i$, whose exceptional set contains a divisor, followed by a small modification $Z^i  \stackrel{\psi_i}{\longleftarrow}X^{i+1}$.
\end{itemize}
Note that the mixed case  can occur  only if  either $X^i$ is not $\q$-factorial or  $\phi_i$  contracts an extremal face of dimension $\geq 2$.
In most treatments this is avoided   by working   with 
$\q$-factorial varieties and  choosing $H$ sufficiently general.

We can almost always choose the initial $X$ to be  nonsingular, but
frequently other considerations constrain  the choice of $H$.  

Our aim is to discuss a significant special case where 
the  $X^i$ are not  $\q$-factorial and we do contract 
 extremal faces of dimension $\geq 2$, but still avoid the  mixed case.  This has several applications, some of which are discussed in Section~2.

\begin{ack} I thank E.~Arvidsson, F.~Bernasconi, J.~Carvajal-Rojas, J.~Lacini,  A.~St{\"a}bler,  D.~Villalobos-Paz, C.~Xu
for  helpful comments and   J.~Witaszek for numerous e-mails about flips.
Partial  financial support    was provided  by  the NSF under grant number
DMS-1901855.
\end{ack}

\section{Relative MMP with scaling of an exceptional divisor}

\begin{defn}[MMP with scaling]\label{mmp.steps.defn}
Let $X,S$ be  Noetherian, normal  schemes and $g:X\to S$  a projective
morphism.  Let $\Theta$ be an  $\r$-divisor on $X$ and 
$H$ an  $\r$-Cartier, $\r$-divisor  on $X$. 
Assume that $K_X+\Theta+r_XH$ is $g$-ample for some $r_X$.

By the  {\it $(X,\Theta)$-MMP with scaling} of  $H$ we mean 
a sequence of normal, projective  schemes  $g_j:X^j\to S$ and birational 
contractions  $\tau_j:X^j\map X^{j+1}$, together with real numbers
$r_X=r_0>r_1\cdots$, that are constructed by 
the following process.

$\bullet$ We start with  $(X^1,\Theta^1, H^1):=(X,\Theta, H)$ and $r_0=r_X$. If $D$ is any divisor on $X$, we let $D^j$ denote its birational transform on $X^j$.

$\bullet$ If $X^j, \Theta^j, H^j$ are already defined, we let  $r_j<r_{j-1}$ be the unique real number for which
 $K_{X^j}+\Theta^j+r_jH^j$ is $g_j$-nef but not $g_j$-ample.
  Then the $j$th step of the MMP is a diagram
$$
\begin{array}{rrcll}
(X^j, \Theta^j) &\stackrel{\phi_j}{\rightarrow} & Z^j & \stackrel{\psi_j}{\leftarrow} &
(X^{j+1}, \Theta^{j+1})\\
 g_j &\searrow &\downarrow & \swarrow &g_{j+1} \\
&& S &&
\end{array}
\eqno{(\ref{mmp.steps.defn}.1)}
$$
where
\begin{enumerate}\setcounter{enumi}{1}
\item  $\phi_j$ is the contraction defined by $K_{X^j}+\Theta^j+r_jH^j$,
\item $\psi_j$ is small, and
\item  $K_{X^{j+1}}+\Theta^{j+1}+(r_j-\epsilon)H^{j+1}$ is  $g_{j+1}$-ample for $0<\epsilon\ll 1$.
\end{enumerate} 
Note that (4) implies that $H^{j+1}$ must be $\r$-Cartier.

In general such a diagram need not exist, but if it does, it is unique and
then $X^{j+1}, \Theta^{j+1}, H^{j+1}$ satisfy the original assumptions.
Thus, as far as the existence of MMP-steps is concerned, we can focus on the 1st step. In this case it is customary to drop the upper indices and write
(\ref{mmp.steps.defn}.1) as
$$
\begin{array}{rrcll}
(X, \Theta) &\stackrel{\phi}{\rightarrow} & Z & \stackrel{\phi^+}{\leftarrow} &
(X^+, \Theta^+)\\
 g &\searrow &\downarrow & \swarrow &g^+ \\
&& S &&
\end{array}
\eqno{(\ref{mmp.steps.defn}.5)}
$$
We say that the MMP {\it terminates} with $g_j:X^j \to S$ if
\begin{enumerate}\setcounter{enumi}{5}
\item either   $K_{X^j}+\Theta^j$ is $g_j$-nef, in which case
$ (X^j,\Theta^j)$ is called a {\it minimal model} of  $(X, \Theta)$,
\item or $\phi_j: X^j{\to} Z^j$  exists and $\dim Z^j<\dim X^j$; then $\phi_j$ is called a {\it Fano contraction}.
\end{enumerate}

{\it Warning  \ref{mmp.steps.defn}.8.}  Our  terminology is slightly different from
\cite{bchm}, where it is  assumed that $X^j/Z^j$ has relative Picard number 1,
and  $r_j=r_{j-1}$ is allowed. In effect, we declare that the  composite of all
\cite{bchm}-steps with the same value of $r$ is a single step for us.
Thus we sometimes contract an extremal face, not just an   extremal ray.

One advantage is that our MMP steps are uniquely determined by
the starting data.  This makes it possible to extend the theory to algebraic spaces \cite{dvp}.

\end{defn}

Theorem~\ref{bir.mmp.thm} is  formulated for  Noetherian base schemes.  We do not prove any new results about the existence of flips, but Theorem~\ref{bir.mmp.thm} says that  if 
the  MMP with scaling exists and terminates, then  its steps are simpler than expected, and the end result is more controlled than expected.

On the other hand, for 3-dimensional schemes, Theorem~\ref{bir.mmp.thm} can be used to conclude that, in some important cases,  the MMP runs and terminates, see Theorem~\ref{mmp.simpler.prop}.

\begin{thm} \label{bir.mmp.thm}
Let $Y$ be a Noetherian, normal  scheme and $g:X\to Y$  a projective, birational morphism with reduced exceptional divisor $E=E_1+\cdots + E_n$. 
 Assume the following
(which are frequently easy to achieve, see Paragraphs~\ref{amp.exc.div.say.0}--\ref{amp.exc.div.say}).
\begin{enumerate}
\item[(i)] $(X,\Theta)$ is dlt and the $E_i$ are $\q$-Cartier.
\item[(ii)]
 $K_X+\Theta\equiv_{g}  E_\Theta$, where $E_\Theta=\sum e_iE_i$.
\item[(iii)] $H=\sum h_iE_i$, where $-H$ is effective and  $\supp H= E=\ex(g)$.
\item[(iv)] $K_X+\Theta+r_XH$ is $g$-ample for some $r_X>0$.
\item[(v)] The $h_i$ are linearly independent over $\q(e_1,\dots, e_n)$.
\end{enumerate}
We run the  $(X,\Theta)$-MMP with scaling of  $H$.
 Assume that we reached the $j$th step  as in (\ref{mmp.steps.defn}.1). Then the following hold.
\begin{enumerate}
\item
 $\ex(\phi_j)\subset \supp(E^j)$ and  
\begin{enumerate}
\item either  $\ex(\phi_j)$ is an irreducible divisor and $X^{j+1}=Z^j$,
\item or $\phi_j$ is  small, and 
there are irreducible components  $E^j_{i_1}, E^j_{i_2}$  of $E^j$   such that  $E^j_{i_1}$ and $ -E^j_{i_2}$ are both $\phi_j$-ample.
\end{enumerate}
\item The $E^{j+1}_i$ are all $\q$-Cartier.
\item  $E_{\Theta}^{j+1}+(r_j-\epsilon)H^{j+1}$ is a $g_{j+1}$-ample $\r$-divisor supported on $\ex(g_{j+1})$ for $0<\epsilon\ll 1$.
\end{enumerate}
\noindent  Furthermore,  if the  MMP terminates with $g_m:X^m\to Y$, then
\begin{enumerate}\setcounter{enumi}{3}
\item    $-E^m_\Theta$ is effective,  $\supp E^m_\Theta=g_m^{-1}\bigl(g_m(\supp E^m_\Theta)\bigr)$, and 
\item if $E^m_\Theta$ is effective and $\supp E_\Theta=E$, then   $X^m=Y$. 
\end{enumerate}
\end{thm}

\noindent  {\it Remark \ref{bir.mmp.thm}.6.} In applications the following are the key points:
\begin{enumerate}
\item[(a)] We avoided the mixed case.
\item[(b)] In the fipping case we have both $\phi$-positive and $\phi$-negative divisors. 
\item[(c)] In (3) we have an explicit, relatively ample, exceptional $\r$-divisor.
\item[(d)] In case (5) we  end with  $X^m=Y$
(not with an unknown,  small modification of $Y$).
\item[(e)]  In case (5) the last MMP step is a divisorial contraction,  giving what \cite{MR3187625}   calls a {\it Koll\'ar component;} no further flips needed.
\end{enumerate}
\medskip

Proof. Assertions (1--3)  concern only one MMP-step, so we may as well 
drop the index $j$ and work with the diagram (\ref{mmp.steps.defn}.5). Thus assume that
  $K_X+\Theta+(r+\epsilon)H$ is $g$-ample,  $K_X+\Theta+rH$ is $g$-nef and 
it determines the  contraction 
$\phi: X\to Z$. 

Let $N_1(X/Z)$ be the relative cone of curves. The $E_i$ give elements of the dual space  $N^1(X/Z)$.  If $C\subset X$ is contracted by $\phi$ then we have a relation  
$$
\tsum h_i(E_i\cdot C)=-r^{-1}(E_\Theta\cdot C).
\eqno{(\ref{bir.mmp.thm}.7)}
$$
 By Lemma~\ref{bir.mmp.thm.lem.1} this shows that the $E_i$ are proportional, as functions on $N_1(X/Z)$. Let $C'$ be another contracted curve; set
$e:=(E_\Theta\cdot C)$ and $e':=(E_\Theta\cdot C')$. Using (\ref{bir.mmp.thm}.7) for $C$ and $C'$, we can eliminate $r$ to get that
$$
\tsum h_i\bigl(e'(E_i\cdot C)-e(E_i\cdot C')\bigr)=0.
\eqno{(\ref{bir.mmp.thm}.8)}
$$
By the linear independence of the $h_i$ this implies that
$e'(E_i\cdot C)=e(E_i\cdot C')$ for every $i$. That is, 
all contracted curves  are proportional, as functions on 
$\langle E_1,\dots, E_n\rangle_{\r}\cong \r^n$.
Informally speaking, as far as the $E_i$ are concerned, $N_1(X/Z)$ behaves as if it were  1-dimensional.

Assume first that $\phi$ contracts some divisor, call it $E_1$. 
Then $(E_1\cdot C)<0$ for some contracted curve  $C\subset E_1$, hence 
$(E_1\cdot C')<0$ for every contracted curve  $C'$. Thus $\ex(\phi_0)=E_1$. We also know that 
$$
\phi_*(E_\Theta+rH)=\tsum_{i>1} (e_i+rh_i) \phi_*(E_i)
$$
 is $\r$-Cartier on $Z$ and $Z/Y$-ample, where
 $r$ is computed by (\ref{bir.mmp.thm}.7). So, 
by Lemma~\ref{bir.mmp.thm.lem.2}, the $\{e_i+rh_i: i>1\}$ are 
linearly independent over $\q$, hence the $\phi_*(E_i)$ are $\q$-Cartier on $Z$
by Lemma~\ref{alex.lem}. Thus   $\phi_*(E_\Theta)=\tsum_{i>1} e_i \phi_*(E_i)$ is $\r$-Cartier, hence
 $X^1=Z$.  This proves (2--3) in the divisorial contraction case.

Otherwise  $\phi$ is small, let $C$ be a contracted curve.  Since $(H\cdot C)>0$, we get that $(E_1\cdot C)<0$ for some  $E_1$. So $C\subset E_1$. 
 By \cite[3.39]{km-book}
$E_\Theta+rH$ is anti-effective and 
$$
g^{-1}\bigl(g(\supp(E_\Theta+rH))\bigr)=\supp(E_\Theta+rH).
\eqno{(\ref{bir.mmp.thm}.9)}
$$
If $E_1$ has coefficient 0 in $E_\Theta+rH$ then let
$C_1\subset E_1$ be any curve contracted by $g$. Then $C_1$ is disjoint from
$\supp(E_\Theta+rH)$ by (\ref{bir.mmp.thm}.9), hence  $(C_1\cdot E_\Theta+rH)=0$. Varying $C_1$ shows that 
 $E_1$ is  contracted by $\phi$, a contradiction. 

Thus
$E_1$ appears in $E_\Theta+rH$ with negative coefficient, contributing a 
positive term to the intersection $\bigl((E_\Theta+rH)\cdot C\bigr)=0$. So there is another divisor $E_2\subset \ex(g)$ such that  $(E_2\cdot C)>0$. This shows (1.b).

Assume next that  the flip  $\phi^+:X^+\to Z$ exists.
Since $\phi^+$ is small, 
$\supp(E^+_\Theta+rH^+)$ contains all $X^+/Y$-exceptional divisors.
In particular, 
 $E^+_\Theta+(r-\epsilon)H^+$ is still  anti-effective 
for $0<\epsilon\ll 1$. 
By definition $E^+_\Theta+(r-\epsilon)H^+$ is  $X^+/Y$-ample 
and its  support   is the whole  $X^+/Y$-exceptional locus. 
Thus we also have  (2--3) in the flipping case.

Finally,  if the  MMP terminates with $g_m:X^m\to Y$ then  
$E^m_\Theta$ is a $g_m$-nef, exceptional $\r$-divisor. 
Thus $-E^m_\Theta$ is effective and 
$\supp E^m_\Theta=g_m^{-1}\bigl(g_m(\supp E^m_\Theta)\bigr)$  by \cite[3.39]{km-book}, proving (4). In case (5) this implies that 
 $\ex(g_m)$ does not contain any divisor, but, 
by (3) it supports a $g_m$-ample divisor. Thus $\dim \ex(g_m)=0$,  hence
   $X^m=Y$. \qed

\begin{lem} \label{bir.mmp.thm.lem.1}
Let $V$ be a $K$-vectorspace with vectors  $v_i\in V$.
Let $L/K$ be a field extension and $h_1,\dots, h_n\in L$ linearly independent over $K$.  Assume that
$$
\tsum_{i=1}^n h_iv_i=\gamma v_0\qtq{for some} \gamma\in L.
$$
Then $\dim_K  \langle  v_1,\dots, v_n\rangle \leq 1$.
\end{lem}

Proof. We may assume that $\dim V=2$.
Choose a basis and write  $v_i=(a_i, b_i)$. Then 
$$
\tsum_{i=1}^n h_ia_i=\gamma a_0\qtq{and}\tsum_{i=1}^n h_ib_i=\gamma b_0.
$$
This gives that
$$
\tsum_{i=1}^n h_i(b_0a_i-a_0b_i)=0.
$$
Since the $h_i$ are  linearly independent over $K$, this implies that
$b_0a_i-a_0b_i=0$ for every $i$. That is
$$
v_i\cdot (b_0, -a_0)^t=0\qtq{for every $i$.} \qed
$$

\begin{lem} \label{bir.mmp.thm.lem.2}
Let $L/K$ be a field extension and $h_0,\dots, h_n\in L$ linearly independent over $K$.    Let $\gamma^{-1}=\tsum_{i=0}^n r_ih_i$ for some $r_i\in K$ with $r_0\neq 0$. 
Then, for any $e_i\in K$, the
$e_1+\gamma h_1,\dots, e_n+ \gamma h_n$ are  linearly independent over $K$. 
\end{lem}

Proof.    Assume that
$\tsum_{i=1}^n s_i (e_i+\gamma h_i)=0$, where $s_i\in K$.
It rearranges to
$$
\tsum_{i=1}^n s_i h_i=-\bigl(\tsum_{i=1}^n s_i e_i\bigr)\cdot  \tsum_{i=0}^n r_ih_i.
$$
If $\tsum_{i=1}^n s_i e_i=0$ then the $s_1,\dots, s_n$ are all zero since the $h_i$ are linearly independent over $K$. 
Otherwise we get  a contradiction since the coefficient of $h_0$ is nonzero. \qed

\begin{lem}\cite[Lem.1.5.1]{MR3380944}\label{alex.lem} Let $X$ be a normal scheme,  $D_i$  $\q$-divisors and $d_1,\dots, d_n\in \r$ linearly independent over $\q$.  Then $\sum d_iD_i$ is $\r$-Cartier iff each $D_i$ is  
$\q$-Cartier. \qed
\end{lem}

\subsection*{Comments on $\q$-factoriality}{\ }
\medskip

Theorem~\ref{bir.mmp.thm} may sound unexpected from the MMP point of view, but it is quite natural if one starts with  the following  conjecture, which is due, in various forms,  to Srinivas  and myself, cf.\ \cite{MR1242007}.

\begin{conj}\label{k-s.conj}
 Let $X$ be a normal variety, $x\in X$ a closed point and $\{D^X_i: i\in I\}$ a finite set of  divisors on $X$. 
Then there is a normal variety $Y$,   a closed point $y\in Y$  and divisors $\{D^Y_i: i\in I\}$ on $Y$ such that the following hold.
\begin{enumerate}
\item The class group of the local ring  $\o_{y,Y}$ is generated by $K_Y$ and the $D^Y_i$.
\item  The completion of $( X, \sum D^X_i)$ at $x$ is isomorphic to the completion of $( Y, \sum D^Y_i)$ at $y$.
\end{enumerate}
\end{conj}

Using \cite[\href{https://stacks.math.columbia.edu/tag/0CAV}{Tag 0CAV}]{stacks-project} one can reformulate  (\ref{k-s.conj}.2) as  a  finite type statement:
\begin{enumerate}\setcounter{enumi}{2}
\item  There are elementary \'etale morphisms
$$(x,X,\tsum D^X_i) \leftarrow (u, U,\tsum D^U_i) \to (y,Y,\tsum D^Y_i).$$
\end{enumerate}

Almost all resolution methods commute with \'etale morphisms, 
thus if we want to prove something about a resolution of $X$, it is likely to be equivalent to a statement about  resolutions of $Y$. 
In particular,  if something holds for the  $\q$-factorial case, it should hold in general.
 This was the reason why I believed that Theorem~\ref{bir.mmp.thm} should work out.

A positive answer to Conjecture~\ref{k-s.conj} (for $I=\emptyset$) is given for isolated complete intersections in  \cite{MR1242007} and for  normal surface singularities in \cite{par-dvs}.

(Note that \cite{par-dvs} uses an even stronger formulation: Every normal, analytic singularity has an algebraization whose class group  is generated by the canonical class. This is, however, not true, since not every normal, analytic singularity has an algebraization.) 

\subsection*{Existence of certain  resolutions}{\ }

\begin{say}[The assumptions \ref{bir.mmp.thm}.i--v]
\label{amp.exc.div.say.0} 
In most applications of Theorem~\ref{bir.mmp.thm} we start with a normal pair 
$(Y, \Delta_Y)$  where $\Delta_Y$  is a boundary,  and want to find $g:X\to Y$ and $\Theta$ that satisfy the conditions (\ref{bir.mmp.thm}.i--v). 

Typically we choose a log resolution $g:X\to (Y, \Delta_Y)$.
That is, $g$ is birational, $X$ is regular, $\Delta_X:=g^{-1}_*\Delta_Y$,
$E=\ex(g)$ and  
$E+\Delta_X$ is a simple normal crossing divisor. 
Then we choose $\Delta_X\leq \Theta\leq E+\Delta_X$; that is, we are free to choose the coefficients of the $E_i$ in $[0,1]$. Then (\ref{bir.mmp.thm}.i) holds and if $K_Y+\Delta_Y$ is $\r$-Cartier then so does (\ref{bir.mmp.thm}.ii). There are also  situations where one can use the theorem to show that numerical equivalence in (\ref{bir.mmp.thm}.ii) implies linear equivalence; see \cite[9.12]{many-p}.

We want $K_X+\Theta +rH$ to be $g$-ample for some $r$, which is easiest to achieve if $H$ is  $g$-ample. Thus we would like $H$ to be $g$-ample and $g$-exceptional for (\ref{bir.mmp.thm}.iii--iv) to hold.
If $X$ is regular  (or at least $\q$-factorial) then we can wiggle
the coefficients of $H$ to achieve   (\ref{bir.mmp.thm}.v).

The existence of a $g$-ample and $g$-exceptional divisor is somewhat subtle, we discuss it next.
\end{say}

\begin{say}[Ample, exceptional divisors]\label{amp.exc.div.say} 
Assume that  we blow up an ideal sheaf $I\subset \o_Y$ to get $\pi_1:Y_1\to Y$.
The constant sections of $\o_Y$ give an isomorphism 
 $\o_{Y_1}(1)\cong \o_{Y_1}(-E_1)$ where $E_1$ is supported on 
$\pi_1^{-1}\supp (\o_Y/I)$. Thus, if $Y$ is normal and
$\supp (\o_Y/I)$  has codimension $\geq 2$, then $E_1$ is 
 $\pi_1$-ample and $\pi_1$-exceptional.
A composite of  such blow-ups also has an ample, exceptional divisor.
Since Hironaka-type resolutions use only such blow-ups, we get the following.
(See \cite{tem} for the most general case and \cite{k-res} for an introduction.)
\medskip

{\it Claim \ref{amp.exc.div.say}.1.} Let $Y$ be a Noetherian, quasi-excellent scheme over a field of  characteristic zero. Then any proper, birational $Y'\to Y$ is dominated by a log resolution $g:X\to Y$ that has a $g$-ample and $g$-exceptional divisor. \qed
\medskip

Resolution of singularities is also known for 3-dimensional excellent schemes \cite{cos-pil-2014},
but in its original form it does not guarantee projectivity in general.
Nonetheless, combining \cite[2.7]{many-p} and 
\cite[Cor.3]{k-wit}  we get the following.
\medskip

{\it Claim \ref{amp.exc.div.say}.2.} 
Let $Y$ be a normal, integral, quasi-excellent scheme  of dimension at most three that is separated and of finite type over an affine, quasi-excellent scheme $S$.   Then any proper, birational $Y'\to Y$ is dominated by a log resolution $g:X\to Y$ that has a $g$-ample and $g$-exceptional divisor. \qed
\end{say}

\section{Applications}

Next we mention some applications. In each case we use Theorem~\ref{bir.mmp.thm} to modify the previous proofs to get more general results.  We give only some hints as to how this is done, we refer to the original papers for definitions and details of proofs.

The first two applications are to  dlt 3-folds. 
In both cases Theorem~\ref{bir.mmp.thm} allows us to run MMP in a way that works in every characteristic and also for bases that are not $\q$-factorial.

\subsection*{Relative MMP for dlt 3-folds}{\ }
\medskip

\begin{thm} \label{mmp.simpler.prop}
 Let $(Y, \Delta)$ be a 3-dimensional, normal, Noetherian, excellent   pair 
such that $K_Y+\Delta$ is $\r$-Cartier and $\Delta$ is a boundary.  Let $g:X\to Y$ be a  log resolution with exceptional divisor $E=\sum E_i$.  Assume that $E$ supports a $g$-ample divisor $H$  (we can then choose its  coefficients  sufficiently general).

Then   the MMP over $Y$, starting with  $(X^0, \Theta^0):=(X, E+g^{-1}_*\Delta)$ with scaling of $H$ runs and  terminates with a minimal model $g_m:(X^m, \Theta^m)\to Y$. Furthermore, 
\begin{enumerate}
\item  each  step  $X^i\map X^{i+1}$ of this MMP is
\begin{enumerate}
\item either a contraction $\phi_i:X^i\to X^{i+1}$, whose exceptional set is an irreducible component of $E^i$,
\item or a flip  $X^i\stackrel{\phi_i}{\longrightarrow}Z^i  \stackrel{\psi_i}{\longleftarrow}(X^i)^+=X^{i+1}$, and  
there are irreducible components  $E^i_{i_1}, E^i_{i_2}$  such that $E^i_{i_1},$ and $ -E^i_{i_2}$ are  both  $\phi_i$-ample,
\end{enumerate}
\item $\ex(g_m)$ supports a $g_m$-ample $\r$-divisor, and
\item if  either  $(Y, \Delta)$ is plt, or   $(Y, \Delta)$ is dlt and $g$ is thrifty \cite[2.79]{kk-singbook}, then $X^m=Y$. 
\end{enumerate}
\end{thm}

Proof.  Assume first that the MMP steps exist and the MMP terminates.
Note that
$$
\begin{array}{rll}
K_X+E+g^{-1}_*\Delta&\simr &g^*(K_Y+\Delta)+\tsum_j\bigl(1+a(E_j, Y, \Delta)\bigr) E_j\\
&\sim_{g, \r} &  \tsum_j\bigl(1+a(E_j, Y, \Delta)\bigr) E_j =:E_\Theta.
\end{array}
$$
 We get from Theorem~\ref{bir.mmp.thm} that (1.a--b) are the possible MMP-steps,  and (2--3)  from
 Theorem~\ref{bir.mmp.thm}.3--5.

For existence and termination, almost everything is covered by  \cite{many-p}.

However, the claim I would like to make is that we are in a much simpler situation, that can be treated with the methods that are already in \cite{sho-3ff, k-etal}.

 The  key point is that everything happens inside $E$.  We can thus understand the whole MMP by looking at the 2-dimensional scheme  $E$. This is easiest for termination, which follows from \cite[Sec.7]{k-etal}.

Contractions for reducible surfaces have been treated in
\cite[Secs.11-12]{k-etal}, see also \cite[Chap.6]{fuj-book} and
\cite{tan-surf}.

The presence of  $E^i_{i_1}, E^i_{i_2}$ means that the flips are rather special; 
 called   {\it 1-complem\-ented flips} in \cite{sho-3ff} and {\it easy flips}    in \cite[Sec.20]{k-etal}.
I believe that the medhods of
\cite{sho-3ff, k-etal}  prove the existence of  1-complemented 3-fold flips if $Y$ is excellent and its closed points have perfect residued fields; but the details have not been written down.

The short note  \cite{jakub-0708} explains how 
 \cite[3.4]{hac-wit1} gives  1-comple\-men\-ted 3-fold flips; see
\cite[3.1 and 4.3]{h-w-p} for stronger results.
\qed

\subsection*{Inversion of adjunction for 3-folds}{\ }
\medskip

Using Theorem~\ref{mmp.simpler.prop} we can remove the $\q$-factoriality assumption from
\cite[Cor.1.5]{hac-wit1}. The characteristic 0 case, in all dimensions, was proved in \cite[17.6]{k-etal},

\begin{cor} Let $(X, S+\Delta)$ be a 3-dimensional, normal, Noetherian, excellent  pair. Assume that $X$ is normal, $S$ is a reduced divisor, $\Delta$ is effective and $K_X+ S+\Delta$ is $\r$-Cartier. Let $\bar S\to S$ denote the normalization. Then 
$(\bar S, \diff_{\bar S}\Delta)$ is klt iff  $(X, S+\Delta)$  is plt  near  $S$. \qed
\end{cor}

This implies that one direction of Reid's classification of terminal singularities using `general elephants'  \cite[p.393]{r-ypg}  works in every characteristic. This could be useful in  extending \cite{2019arXiv191208779A} to characteristics $\geq 5$. 

\begin{cor} Let $(X, S)$ be a 3-dimensional pair. Assume that $X$ is normal, $K_X+ S$ is Cartier, $X$ and $S$ have only isolated singularities, and  the normalization
$\bar S$ of $S$ has canonical singularities.  Then $X$  has terminal singularities  in a neighborhood of $S$. \qed
\end{cor}

\subsection*{Divisor class group of dlt singularities}{\ }
\medskip

The divisor class group of a rational surface singularity is finite
by \cite{lip-rs}, and \cite{MR492263}  plus an easy argument shows that the divisor class group of a rational 3-dimensional  singularity is finitely generated. Thus the divisor class group of a  3-dimensional  dlt singularity is finitely generated  in characteristic  $\geq 7$, using  \cite[Cor.1.3]{2020arXiv200603571A}.   Theorem~\ref{mmp.simpler.prop} leads---via \cite[B.14]{k-crst}---to the following weaker result, which is, however, optimal  in characteristics $2,3,5$; see \cite{cr-st} for an application.

\begin{prop}\cite[B.1]{k-crst} \label{dlt.loc.pic.prop}
Let $(y,Y, \Delta)$ be a  3-dimensional,  Noetherian, excellent, dlt singularity with residue characteristic $p>0$.  Then the prime-to-$p$ parts of  $\cl(Y), \cl(Y^{\rm h})$ and of   $\cl(\hat Y)$  are
 finitely generated, where $Y^{\rm h}$ denotes the henselisation and 
$\hat Y$ the completion. \qed
\end{prop}

It seems reasonable to conjecture that the same holds in all dimensions, see 
\cite[B.6]{k-crst}.

\subsection*{Grauert-Riemenschneider vanishing}{\ }
\medskip

One can prove a variant of the Grauert-Riemenschneider (abbreviated as G-R) vanishing theorem \cite{Gra-Rie70b} by following the steps of the MMP.

\begin{defn}[G-R vanishing]\label{g-r.prop.defn}
Let $(Y, \Delta_Y)$ be a pair, $Y$ normal, $\Delta_Y$ a boundary (that is, all coefficients are in $[0,1]$) and $g:X\to Y$ a proper, birational morphism with $X$ normal. 
For an $\r$-divisor  $F$ on $X$ let $\ex(F)$ denote its $g$-exceptional part.
Assume that $Y$ has a dualizing complex.
We say that  {\it  G-R vanishing} holds for   $g:X\to (Y,\Delta_Y)$ if the following is  satisfied.

Let  $D$ be a  $\z$-divisor  and $\Delta_X$ an  effective  $\r$-divisor on $X$. Assume that  
\begin{enumerate}
\item    $D\sim_{g,\r} K_{X}+\Delta_X$, and
\item      $g_*\Delta_X\leq \Delta_Y$,  $\rdown{\ex(\Delta_X)}=0$.
\end{enumerate}
Then  $R^ig_*\o_{X}(D)=0$ for $i>0$.

We say that  {\it  G-R vanishing} holds over   $(Y, \Delta_Y)$ if
 G-R vanishing holds for every log resolution $g:(X, E+g^{-1}_*\Delta_Y)\to (Y, \Delta_Y)$.

 By an elementary computation,
if $X$ is regular, $W\subset X$ is regular  and G-R vanishing holds  for $X\to Y$ then it also holds for  the blow-up $B_WX\to Y$. This implies that if  G-R vanishing holds  for one log resolution of $(Y,\Delta_Y)$, then it holds for every log resolution; see \cite[Sec.1.3]{k-ber}.

If $Y$ is essentially of finite type over a field of characteristic 0, then
 G-R vanishing is a special case of the general Kodaira-type vanishing theorems; see \cite[2.68]{km-book}.  

 G-R vanishing also holds over 2-dimensional, excellent schemes by
\cite{lip-rs}; see \cite[10.4]{kk-singbook}.
In particular, if $Y$ is any normal, excellent scheme,  then
the support of $R^ig_*\o_{X}(D)=0$ has codimension $\geq 3$ for $i>0$.

However, G-R vanishing fails for 3-folds in every positive characteristic, as shown by cones over surfaces for which Kodaira's vanishing fails. Thus the following may be the type of  G-R vanishing result that one can hope for.
\end{defn}

\begin{thm}\cite{k-ber} \label{dlt.chr7.gr.cor}
Let $Y$ be a 3-dimensional, excellent,  dlt pair with a dualizing complex.  
Assume that closed points of $Y$ have  perfect residue fields  of    characteristic  $\neq 2,3,5$.  Then G-R vanishing holds over  $Y$.
\end{thm}

Proof. Let $(Y, \Delta_Y)$ be a 3-dimensional,  dlt pair, and $g:X\to Y$ a log
resolution. With $D$ as in Definition~\ref{g-r.prop.defn} we need to show that
 $R^jg_*\o_{X}(D)=0$ for $j>0$.  Let $g_i:X^i\to Y$ be the MMP steps as in  Theorem~\ref{mmp.simpler.prop}.  The natural idea would be to show that
the sheaves  $R^j(g_i)_*\o_{X^i}(D^i)$ are independent of $i$. At the end then we have an isomorphism $g_m:X^m\cong Y$, hence $R^j(g_m)_*\o_{X^m}(D^m)=0$ for $j>0$.

A technical problem is that we seem to need various rationality properties of the singularities of the $X^i$. Therefore,   we show instead that, if 
G-R vanishing holds {\em over} $X^i$ 
and  $X^i$ satisfies (\ref{g-r.impl.other.thm}.1--3), then G-R vanishing also holds {\em over} $X^{i+1}$.  Then Theorem~\ref{g-r.impl.other.thm} gives that 
 $X^{i+1}$  also satisfies  (\ref{g-r.impl.other.thm}.1--3), and the induction can go ahead.

For divisorial contractions $X^i\to X^{i+1}$ with exceptional divisor $S$ this is straightfoward, the method
of \cite[Sec.3]{h-w-1} shows that if Kodaira vanishing holds  for $S$ then
G-R vanishing holds for $X^i\to X^{i+1}$. This is where the 
$\chr \neq 2,3,5$ assumption is used: Kodaira vanishing can fail for
del~Pezzo surfaces if $\chr =2,3, 5$; see \cite{2020arXiv200603571A}. 

For flips  $X^i\to Z^i\leftarrow  X^{i+1}$ the argument works in any characteristic. First we show as above that G-R vanishing holds over $Z^i$.
Going to $X^{i+1}$ is a spectral sqeuence argument involving
$\psi_i:X^{i+1}\to Z^i$. For 3-folds the only nontrivial term is
$R^1(\psi_i)_*\o_{X^{i+1}}(D^{i+1})$, and no unexpected cancellations  occur; see \cite[Lem.21]{k-ber}. \qed
\medskip

From G-R vanishing  one can derive  various rationality properties for all excellent   dlt pairs.  This can be done by following the method of 2 spectral sequences  as in  \cite{k-dep} or \cite[7.27]{kk-singbook}; see
\cite{k-ber} for an improved version.

\begin{thm}\cite{k-ber}\label{g-r.impl.other.thm}
 Let $(X, \Delta)$ be an excellent  dlt pair such that  G-R vanishing  and resolution of singularities hold over   $(X, \Delta)$. Then
\begin{enumerate}
\item $X$ has rational singularities.
\item Every irreducible component of $\rdown{\Delta}$ is normal and has rational singularities.
\item Let $D$ be a   $\z$-divisor on $X$ such that
$D+ \Delta_D$ is $\r$-Cartier for some $0\leq \Delta_D\leq \Delta$. Then 
$\o_X(D)$ is CM. \qed
\end{enumerate}
\end{thm}

See \cite[12]{k-ber} for the precise resolution assumptions needed.
The conclusions are well known in characteristic 0, see   
\cite[5.25]{km-book}, \cite[Sec.3.13]{fuj-book}  and \cite[7.27]{kk-singbook}.
 For 3-dimensional  dlt varieties in $\chr \geq 7$,  the first claim was  proved  in \cite{h-w-1, 2020arXiv200603571A}.

\medskip

The next two applications are in characteristic 0.

\subsection*{Dual complex of a resolution}{\ }
\medskip

Our results can be used to remove the $\q$-factoriality assumption from 
\cite[Thm.1.3]{dkx}. We refer to \cite{dkx} for the definition of a dual complex and the notion of collapsing of a regular cell complex.
We start with the weaker form, Corollary~\ref{collapse.cor}, and then state and outline the proof of the stronger version, Theorem~\ref{collapse.thm}.

\begin{cor}\label{collapse.cor}  Let $(Y, \Delta)$ be a dlt variety over field of characteristic 0 and
$g:X\to Y$ a thrifty log resolution 
whose exceptional set supports a $g$-ample divisor.
For a closed point $y\in Y$ let $E_y\subset g^{-1}(y)$ denote the  divisorial part. Then
 ${\mathcal D}(E_y)$ is collapsible to a point (or it is empty). 
\end{cor}

\begin{thm}\label{collapse.thm} Let $(Y, \Delta)$ be a dlt variety over field of characteristic 0 and
$g:X\to Y$ a projective, birational morphism with exceptional set $E=\cup_i E_i$. For $y\in Y$ let $E_y\subset g^{-1}(y)$ denote the  divisorial part.
Assume that
\begin{enumerate}
\item  $(X, E+g^{-1}_*\Delta)$ is dlt and the  $E_i$ are $\q$-Cartier.
\item $a(E_i, Y, \Delta)>-1$  for every $i$.
\item $E$ supports a $g$-ample divisor.
\end{enumerate}
 Then
 ${\mathcal D}(E_y)$ is collapsible to a point (or it is empty). 
\end{thm}

Proof. Fix $y\in Y$. We may assume that $(y, Y)$ is local and, after passing to an elementary \'etale neighborhood
(cf.\ \cite[\href{https://stacks.math.columbia.edu/tag/02LD}{Tag 02LD}]{stacks-project})  of $y\in Y$,
we may also assume that
$g^{-1}(y)\cap E_i$ is connected  for every irreducible exceptional divisor $E_i$ (cf.\ \cite[\href{https://stacks.math.columbia.edu/tag/04HF}{Tag 04HF}]{stacks-project}).

Let us now run the $(X, E+g^{-1}_*\Delta)$-MMP with scaling of a $g$-ample $\r$-divisor $H$ that is supported on $E$ and has sufficiently general  coefficients.  Theorem~\ref{bir.mmp.thm} applies, as we observed during the proof of Theorem~\ref{mmp.simpler.prop}.

Note that ${\mathcal D}(E_y)\subset {\mathcal D}(E)$ is a {\it full} subcomplex
(that is, a simplex is in ${\mathcal D}(E_y)$ iff all of its vertices are), hence an elementary  collapse of ${\mathcal D}(E)$ induces an
elementary  collapse (or an isomorphism) on ${\mathcal D}(E_y)$.
Thus it is enough to show that ${\mathcal D}(E)$ is collapsible to a point (or it is empty). 

We claim that each MMP-step as in Theorem~\ref{bir.mmp.thm}
induces either a  collapse or an isomorphism of ${\mathcal D}(E)$.

By \cite[Thm.19]{dkx} we get an elementary collapse (or an isomorphism) if
there is a divisor $E^j_i\subset E^j$ that has positive intersection with the $\phi_j$-contracted curves. This takes care of flips by Theorem~\ref{bir.mmp.thm}.1.b and most divisorial contractions.

It remains to deal with the case when we contract $E^j_\ell\subset E^j$
and every other 
$E^j_i$ has 0 intersection number with the contracted curves.
Thus $E^j_i\cap E^j_\ell$ is either empty or contains
$g_j^{-1}(y)\cap E^j_\ell$.  Thus the link of $E^j_\ell$ in 
${\mathcal D}(E^j)$ is a simplex and removing it is a sequence of
elementary collapses.  \qed

\subsection*{Dlt modifications of algebraic spaces}{\ }
\medskip

By \cite{MR2955764}, a normal, quasi-projective pair 
$(X, \Delta)$ (over a field of characteristic 0) has both  dlt and lc  modifications if  $K_X+\Delta$ is $\r$-Cartier. 
(See \cite[Sec.1.4]{kk-singbook} for the definitions.)
The lc  modification is unique and commutes with \'etale base change, hence local lc  modifications automatically glue to give the same conclusion if $X$ is an algebraic space.

However, dlt  modifications are rarely unique, thus it was not obvious that they exist when the base is not  quasi-projective. 
\cite{dvp} observed that Theorem~\ref{bir.mmp.thm} gives enough uniqueness to allow for gluing. This is not hard when $X$ is a scheme, but needs careful considerations  to work for algebraic spaces.

\begin{thm}[Villalobos-Paz] Let $X$ be a normal algebraic space of finite type over a field of characteristic 0, and $\Delta$ a boundary $\r$-divisor on $X$. Assume that $K_X+\Delta$ is $\r$-Cartier. Then $(X, \Delta)$ has a modification
$g: (X^{\rm dlt}, \Delta^{\rm dlt})\to (X, \Delta)$ such that
\begin{enumerate}
\item $(X^{\rm dlt}, \Delta^{\rm dlt}) $ is dlt,
\item $K_{X^{\rm dlt}}+ \Delta^{\rm dlt}$ is $g$-nef,
\item $g_*\Delta^{\rm dlt}=\Delta$,  and
\item $g$ is projective.
\end{enumerate}
$X^{\rm dlt}$ is not unique, and we can choose
\begin{enumerate}\setcounter{enumi}{4}
\item   either $X^{\rm dlt}$ to be $\q$-factorial, or
$\ex(g)$ to support a $g$-ample divisor.
\end{enumerate}
\end{thm}


\newcommand{\etalchar}[1]{$^{#1}$}
\def\cprime{$'$} \def\cprime{$'$} \def\cprime{$'$} \def\cprime{$'$}
  \def\cprime{$'$} \def\cprime{$'$} \def\dbar{\leavevmode\hbox to
  0pt{\hskip.2ex \accent"16\hss}d} \def\cprime{$'$} \def\cprime{$'$}
  \def\polhk#1{\setbox0=\hbox{#1}{\ooalign{\hidewidth
  \lower1.5ex\hbox{`}\hidewidth\crcr\unhbox0}}} \def\cprime{$'$}
  \def\cprime{$'$} \def\cprime{$'$} \def\cprime{$'$}
  \def\polhk#1{\setbox0=\hbox{#1}{\ooalign{\hidewidth
  \lower1.5ex\hbox{`}\hidewidth\crcr\unhbox0}}} \def\cdprime{$''$}
  \def\cprime{$'$} \def\cprime{$'$} \def\cprime{$'$} \def\cprime{$'$}
\providecommand{\bysame}{\leavevmode\hbox to3em{\hrulefill}\thinspace}
\providecommand{\MR}{\relax\ifhmode\unskip\space\fi MR }
\providecommand{\MRhref}[2]{%
  \href{http://www.ams.org/mathscinet-getitem?mr=#1}{#2}
}
\providecommand{\href}[2]{#2}

\bigskip

  Princeton University, Princeton NJ 08544-1000, \

\email{kollar@math.princeton.edu}

\end{document}